\documentclass[12pt]{article}
\usepackage{amsmath,amsfonts,amssymb,amscd}
\usepackage[enableskew]{youngtab}
\usepackage{theorem}
\usepackage[usenames]{color}
\nonstopmode

\newtheorem{theorem}{Theorem}[section]

\newtheorem{lemma}[theorem]{Lemma}
\newtheorem{definition}[theorem]{Definition}

\theorembodyfont{\rmfamily}
\newtheorem{example}[theorem]{Example}
\newtheorem{remark}[theorem]{Remark}

\begin{document}
$\,$\vspace{10mm}

\begin{center}
{\textsf {\huge Relationships Between Two Approaches:}}
\vspace{3mm}\\
{\textsf {\huge Rigged Configurations and 10-Eliminations}}
\vspace{15mm}\\
{\textsf{\LARGE  ${}^{\mbox{\small a}}$Anatol N. Kirillov and
${}^{\mbox{\small b}}$Reiho Sakamoto}}
\vspace{20mm}\\
{\textsf {${}^{\mbox{\small{a}}}$Research Institute for Mathematical Sciences,}}
\vspace{-1mm}\\
{\textsf {Kyoto University, Sakyo-ku,}}
\vspace{-1mm}\\
{\textsf {Kyoto, 606-8502, Japan}}
\vspace{-1mm}\\
{\textsf {kirillov@kurims.kyoto-u.ac.jp}}
\vspace{3mm}\\
{\textsf {${}^{\mbox{\small{b}}}$Department of Physics, Graduate School of Science,}}
\vspace{-1mm}\\
{\textsf {University of Tokyo, Hongo, Bunkyo-ku,}}
\vspace{-1mm}\\
{\textsf {Tokyo, 113-0033, Japan}}
\vspace{-1mm}\\
{\textsf {reiho@spin.phys.s.u-tokyo.ac.jp}}
\vspace{40mm}
\end{center}

\begin{abstract}
\noindent
There are two distinct approaches to
the study of initial value problem of the periodic box-ball systems.
One way is the rigged configuration approach
due to Kuniba--Takagi--Takenouchi and another way
is the 10-elimination approach due to
Mada--Idzumi--Tokihiro.
In this paper, we describe precisely interrelations
between these two approaches.
\bigskip\\
Mathematics Subject Classification (2000)
17B37, 37K15, 05E15.\\
Key words and phrases: crystal basis, periodic box-ball system,
combinatorics.
\end{abstract}

\pagebreak

\section{Introduction}
The main purpose of the present paper is to
compare two approaches to the study of dynamics
of the periodic box-ball systems (PBBS for short)
of type $A^{(1)}_1$ due to Kuniba--Takagi--Takenouchi
\cite{KTT} (KTT for short) and Mada--Idzumi--Tokihiro
\cite{MIT} (MIT for short). Our main result states that
$$\mathrm{KTT}\approx\mathrm{MIT}.$$
The approach developed in \cite{KTT} is based on
the theory of the rigged configurations (RC for short ---
for details concerning 
the RC-bijection, see e.g. \cite{KR,S0}),
whereas the approach developed by \cite{MIT}
is based on the 10-elimination procedure.
To be more specific, our main result
describes precisely interrelations
between the 10-elimination algorithm
and the RC-bijection in the case under consideration.

{\it Brief history.}
The box-ball systems (BBS for short)
have been introduced by Takahashi--Satsuma
in 1990 \cite{TS,T} during their study
of cellular automaton and attempts to construct examples of those
which have a solitonical nature.
The periodic version of BBS, PBBS, was then introduced
in \cite{YT,YYT}.
{}From that time the BBS were extensively studied
since many deep and unexpected connections
with different branches of mathematics
and mathematical physics were discovered.
Among those are connections with the theory of
crystal base \cite{HHIKTT,FOY},
combinatorics \cite{TTS,Fuk,KS1},
the Riemann theta functions \cite{KS2,KS3},
tropical algebraic geometry \cite{IT1,IT2},
and also with the theory of discrete KP
and Toda type integrable systems \cite{TTMS,Y,KSY}.

The organization of the present paper is as follows.
In Section \ref{sec:10elimination},
we briefly recall the 10-elimination procedure.
In Section \ref{sec:RC-KTT}, we recall the
rigged configuration and initial value problem
of the PBBS in terms of the RC-bijection.
In Section \ref{sec:main}, we give a precise
interpretation of the 10-elimination in terms
of the RC-approach (Theorem \ref{th:main}).

\section{10-elimination}\label{sec:10elimination}
A state of the PBBS is given by a sequence of
integers 0 and 1 on a circle.
We usually cut the circle at a suitable position
to be determined below and treat it as a sequence on a straight line
which we call a path.
In this paper, we always assume that within each path,
the number of letters 0 is equal to or bigger than
the number of letters 1.
We will call such paths as positive weight paths.
Due to the second paragraph of Section 3.3 of \cite{KTT},
this assumption does not result in loss of generality.
Originally, a path is identified with a sequence of
capacity one boxes, and 0 stands for a vacant box
and 1 stands for a ball within the box.
Given a path $p$, we introduce 10-elimination procedure
and notion of 0-solitons \cite{MIT}.
Input of the 10-elimination procedure is a path $p$ and
output is a finite sequence of paths $E^0(p):=p$, $E^1(p)$,
$E^2(p)$, $\cdots$ defined recursively as follows.

Suppose that we have constructed $E^{k-1}(p)$.
We give coordinate 1, 2, 3, $\cdots$,
to each letter in $E^{k-1}(p)$ from left to right.
The 10-pair is a neighboring pair of 10 where
1 at location $j$ and 0 at location $j+1$ for some $j$.
Then $E^{k}(p)$ is obtained by erasing all
10-pairs of $E^{k-1}(p)$.
If $E^{k}(p)$ does not contain letter 1 for some $k$,
stop the procedure.
Let the number of 10-pairs in $E^{k}(p)$ be $e_k$.
By definition, we have $e_{k-1}\geq e_k$.
Then we denote distinct lengths of rows
of the transposed diagram
${}^t(e_1,e_2,\cdots)$ by $L_j$.
We assume that $L_1>L_2>L_3>\cdots >L_s$
and denote the multiplicity of $L_j$ by
$m_{L_j}$, i.e.,
\begin{align}
{}^t(e_1,e_2,\cdots)=
(L_1^{m_{L_1}},L_2^{m_{L_2}},\cdots,L_s^{m_{L_s}}).
\end{align}

Position of the cut of the original circle
is determined by the following condition.
Note that to each letter in $E^{k}(p)$,
we can specify the original position of $p$
from which the letter is originated.
For a 10-pair in some $E^{k}(p)$,
we join the corresponding 1 and 0 in $p$
by arc in this direction.
Since the path is a positive weight path,
we can always choose a suitable cyclic shift so that no such
arc cross the left or the right end of $p$.
In the sequel we assume that the path $p$
is cut with this property.

Given a path $p$, we draw all possible arcs according to the
procedure in the last paragraph.
Then one step time evolution of the path $p$
is obtained by replacing all connected 1 and 0
by 0 and 1, respectively,
and leaving non-connected letters unchanged.
We denote the resulting time evolved path by $T_\infty (p)$.

Finally, we introduce a notion of 0-solitons.
As we have seen, 10-pairs of $E^{k-1}(p)$ are
erased in $E^k(p)$.
In the following diagram,
$m(>0)$ 10-pairs between $X$ and $Y$ are erased:
\begin{align*}
E^{k-1}(p)=\cdots X(10)^mY\cdots
\xrightarrow{\mbox{ 10-elimination }}
E^{k}(p)=\cdots XY\cdots .
\end{align*}
Under this setting,
\begin{enumerate}
\item[(a)] if $XY=11,01,00$, then there are
$m$ 0-solitons at position $X$ of $E^k(p)$,
\item[(b)] if $XY=10$, then there are $(m-1)$ 0-solitons
at position $X$ of $E^k(p)$.
\end{enumerate}
We can grasp the meaning of the 0-soliton if
we consider the inverse procedure to get
$E^{k-1}(p)$ from $E^k(p)$:
0-solitons give information about how many extra
10-pairs should be inserted between $XY$
to get $E^{k-1}(p)$.
Here, notice that if $XY=10$, it is guaranteed
that there is at least
one 0-soliton between $XY$ in $E^{k-1}(p)$
so that we have no need to specify it.
Based on this observation,
we see that there are precisely $m_{L_j}$
0-solitons at $E^{L_j}(p)$.
Denote the positions of 0-solitons of
$E^{L_j}(p)$ by
\begin{align}
x_1^{(j)}\leq x_2^{(j)}\leq\cdots\leq x_{m_{L_j}}^{(j)}.
\end{align}
As we will see, the data $L_j$ combined with
$x_i^{(j)}$ provide sufficient information
to solve initial value problem of the PBBS.

\begin{example}\label{ex:MIT}
Let us take the length 32 path
$$p=00111011100100011110001101000000.$$
Then 10-elimination procedure goes as follows.
\begin{center}
\setlength{\tabcolsep}{2pt}
\begin{tabular}{lcccccccccccccccccccccccccccccccc}
$E^0(p)=$&0&0&1&1&1&0&1&1&1&0&0&1&0&0&0&1&1&1&1&0&0&0&1&1&0&1&0&0&0&0&0&0,\\
$E^1(p)=$&0&0&1&1& & &1&1& & &0& & &0&0&1&1&1& & &0&0&1& & & & &0&0&0&0&0,\\
$E^2(p)=$&0&0&1&1& & &1& & & & & & &0&0&1&1& & & & &0& & & & & & &0&0&0&0,\\
$E^3(p)=$&0&0&1&1& & & & & & & & & & &0&1& & & & & & & & & & & & &0&0&0&0,\\
$E^4(p)=$&0&0&1& & & & & & & & & & & & & & & & & & & & & & & & & & &0&0&0,\\
$E^5(p)=$&0&0& & & & & & & & & & & & & & & & & & & & & & & & & & & & &0&0.
\end{tabular}
\end{center}
{}From these data, we obtain the following data:
\begin{align}
&(e_0,e_1,e_2,e_3,e_4)=(6,3,2,2,1),\\
&{}^t(e_0,e_1,e_2,e_3,e_4)=
(L_1^{m_{L_1}},L_2^{m_{L_2}},L_3^{m_{L_3}},L_4^{m_{L_4}})=
(5^1,4^1,2^1,1^3),\label{MITdata1}\\
&\{x^{(1)}_1,x^{(2)}_1,x^{(3)}_1,x^{(4)}_1,x^{(4)}_2,x^{(4)}_3\}=
\{2,3,10,4,7,15\}.\label{MITdata2}
\end{align}
We draw arcs on $p$ as follows:
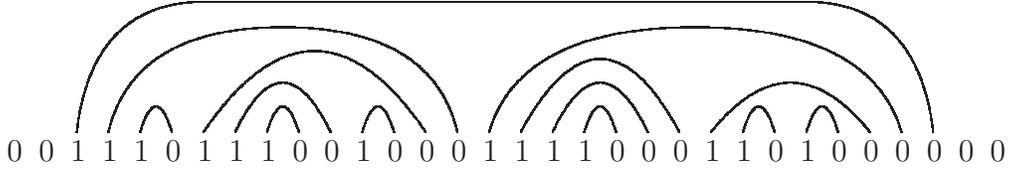
\begin{center}
\unitlength 12pt
\begin{picture}(33,5)(0.5,0.3)
\put(1,0){0}
\put(2,0){0}
\put(3,0){1}
\qbezier(3.2,0.9)(3.7,5)(7.2,5)
\put(4,0){1}
\qbezier(4.2,0.9)(5,4.2)(9.7,4.2)
\put(5,0){1}
\qbezier(5.2,0.9)(5.7,2.5)(6.2,0.9)%
\put(6,0){0}
\put(7,0){1}
\qbezier(7.2,0.9)(10.7,6)(14.2,0.9)
\put(8,0){1}
\qbezier(8.2,0.9)(9.7,4)(11.2,0.9)
\put(9,0){1}
\qbezier(9.2,0.9)(9.7,2.5)(10.2,0.9)%
\put(10,0){0}
\put(11,0){0}
\put(12,0){1}
\qbezier(12.2,0.9)(12.7,2.5)(13.2,0.9)%
\put(13,0){0}
\put(14,0){0}
\put(15,0){0}
\qbezier(15.2,0.9)(14.4,4.2)(9.7,4.2)
\put(16,0){1}
\qbezier(16.2,0.9)(17.2,4.2)(22.7,4.2)
\put(17,0){1}
\qbezier(17.2,0.9)(19.7,5.5)(22.2,0.9)
\put(18,0){1}
\qbezier(18.2,0.9)(19.7,4)(21.2,0.9)
\put(19,0){1}
\qbezier(19.2,0.9)(19.7,2.5)(20.2,0.9)%
\put(20,0){0}
\put(21,0){0}
\put(22,0){0}
\put(23,0){1}
\qbezier(23.2,0.9)(25.7,4)(28.2,0.9)
\put(24,0){1}
\qbezier(24.2,0.9)(24.7,2.5)(25.2,0.9)%
\put(25,0){0}
\put(26,0){1}
\qbezier(26.2,0.9)(26.7,2.5)(27.2,0.9)%
\put(27,0){0}
\put(28,0){0}
\qbezier(29.2,0.9)(28.2,4.2)(22.7,4.2)
\put(29,0){0}
\qbezier(30.2,0.9)(29.7,5)(26.2,5)
\qbezier(7.2,5)(16.7,5)(26.2,5)
\put(30,0){0}
\put(31,0){0}
\put(32,0){0}
\end{picture}
\end{center}
We reverse all connected 10-pairs and obtain
$$T_\infty (p)=00000100011011100001110010111100$$
as the one step time evolution.
\hfill$\square$
\end{example}

\section{Rigged configurations and the KTT theorem}
\label{sec:RC-KTT}
\subsection{Rigged configuration}
Let $X=\{x_1<x_2<\cdots <x_N \}$
be an ordered set, and $\alpha$ be a composition of size $N$.
A standard tabloid of shape $\alpha$ is a filling of the
shape $\alpha$ by elements of the ordered set $X$
such that the elements in each row are strictly increasing.
Let $\nu$ be a partition.
Then $m_i(\nu)$ is the number of occurrences of
$i$ in $\nu$ and $\ell (\nu)$ is the length of $\nu$.
\begin{definition}
A path of type  $B_1^{\otimes L}$
$(=$ path of length $L$ for short$)$
is a sequence $p=a_1a_2\cdots a_L$
where $a_i\in\{0,1\}$ for all $1\leq i\leq L$.
To each path $p=a_1\cdots a_L$ we associate a two row
standard Young tabloid according to the following rule:
we put $i$ to the first row if $a_i=0$, and to
the second row otherwise.
\hfill$\square$
\end{definition}

Our first goal is to remind a construction of the rigged
configuration bijection in the special case of $A^{(1)}_1$.
An input of the RC-bijection is a standard tabloid $T$
of shape $(\lambda_1,\lambda_2)$.
The output of the RC-bijection is a rigged partition,
i.e., a partition $\nu =(\nu_1,\nu_2,\cdots)$ of size $\lambda_2$
together with a collection of integer numbers
\begin{align}
\{J_{\alpha,i}\, |\,
1\leq i\leq\nu_1 \mbox{ s.t. } m_i(\nu)\neq 0,
1\leq\alpha\leq m_i(\nu)\}
\end{align}
such that
\begin{align}
-i\leq J_{\alpha,i}\leq P_i(\nu)
:= L-2Q_i(\nu),
\end{align}
and certain additional restrictions that are not
important for our considerations.
Here $Q_i(\nu):=\sum_a\min (i,\nu_a)$ represents
number of boxes contained in the left $i$ columns
of diagram $\nu$.
The integers $J_{\alpha ,i}$ are called the {\bf riggings}
and $P_i(\nu)$ are called the {\bf vacancy numbers}.

Let us briefly remind a construction of a map from the set of
two row standard tabloids to the set of the rigged partitions.
The construction runs as follows.
Let
$\newcommand{\bone}{b_1}
\newcommand{\blambda}{b_{\lambda_2}}
\Yboxdim 16pt
\Yvcentermath1\young(\bone\cdots\cdots\blambda)$
be the second row of a given two row standard tabloid.
The first step of our construction is to consider
tabloid $\newcommand{\bone}{b_1}
\Yvcentermath1\young(\bone)$
and define the corresponding rigged partition to be
$(\Yvcentermath1\yng(1)\, J_{1,1})$ where $J_{1,1}=b_1-2$.
The next step is to consider tabloid
$\Yvcentermath1
\newcommand{\bone}{b_1}
\newcommand{\btwo}{b_2}
\young(\bone\btwo)$
and define the corresponding rigged partition to be
$(\Yvcentermath1\yng(2)\, J_{1,2})$ where $J_{1,2}=b_2-4$
if $b_2-b_1=1$, and
({\unitlength 12pt
\begin{picture}(2.8,2)
\put(0,0){$\Yvcentermath1\yng(1,1)$}
\put(1.3,0.5){$J_{2,1}$}
\put(1.3,-0.7){$J_{1,1}$}
\end{picture}})
where $J_{1,1}=b_1-2$
and $J_{2,1}=b_2-4$ if $b_2-b_1>1$.
We emphasize that in the case $b_2-b_1=1,$ the 
first row of the configuration
$\nu:=\Yvcentermath1\yng(1,1)$ is singular, i.e., it 
contains the rigging with the maximal possible value.  

We proceed further by induction.
Assume that tabloid
$\newcommand{\bone}{b_1}
\newcommand{\bk}{b_{k}}
\Yboxdim 16pt
\Yvcentermath1\young(\bone\cdots\cdots\bk)$
$(1\leq k<\lambda_2)$,
corresponds to the rigged partition
\begin{align}
\{\tilde{\nu},\tilde{J}_{\alpha ,i}\, |\,
1\leq \alpha\leq m_i(\tilde{\nu})\},
\qquad
\, |\,\tilde{\nu}\, |\,=k.
\end{align}
The next step is to describe the rigged partition
$\{\nu,J_{\alpha,i}\}$
that corresponds to tabloid
$\newcommand{\bone}{b_1}
\newcommand{\bk}{b_{k}}
\newcommand{\bkone}{b_{k+1}}
\Yboxdim 21pt
\Yvcentermath1\young(\bone\cdots\cdots\bk\bkone)$.
If $b_{k+1}-b_k>1$, or equivalently $a_{b_{k+1}-1}=0$,
then $\nu =(\tilde{\nu},1)$ and
$J_{\alpha ,i}=\tilde{J}_{\alpha ,i}$ if $i>1$,
or $i=1$ and $\alpha\leq m_1(\tilde{\nu})$,
and
$J_{m_1(\nu),1}=b_{k+1}-2(\ell (\tilde{\nu})+1)
=\tilde{J}_{m_1(\tilde{\nu})}+b_{k+1}-b_k-2$.
In the case $b_{k+1}-b_k=1$,
there should exist non empty singular strings,
i.e., rows whose riggings are equal to the corresponding
vacancy numbers.
Take a singular string, say $\tilde{\nu}_a$,
of the longest length.
Then the corresponding rigged partition
is defined to be $\{\nu ,J_{\alpha ,i}\}$,
where $\nu_i=\tilde{\nu_i}$ if $i\neq a$, and
$\nu_a=\tilde{\nu_a}+1$.
Moreover, $J_{\alpha ,i}=\tilde{J}_{\alpha ,i}$
if $i\neq a,$ $J_{m_{\nu_a}(\nu),\nu_a}
=b_{k+1}-2\sum_j\min (\nu_a,\nu_j)$.

\begin{example}\label{ex:MIT2}
Let us consider the path $p$ treated in Example \ref{ex:MIT}.
The corresponding tabloid is
$$
\newcommand{\onezero}{10}
\newcommand{\oneone}{11}
\newcommand{\onetwo}{12}
\newcommand{\onethree}{13}
\newcommand{\onefour}{14}
\newcommand{\onefive}{15}
\newcommand{\onesix}{16}
\newcommand{\oneseven}{17}
\newcommand{\oneeight}{18}
\newcommand{\onenine}{19}
\newcommand{\twozero}{20}
\newcommand{\twoone}{21}
\newcommand{\twotwo}{22}
\newcommand{\twothree}{23}
\newcommand{\twofour}{24}
\newcommand{\twofive}{25}
\newcommand{\twosix}{26}
\newcommand{\twoseven}{27}
\newcommand{\twoeight}{28}
\newcommand{\twonine}{29}
\newcommand{\threezero}{30}
\newcommand{\threeone}{31}
\newcommand{\threetwo}{32}
\young(126\onezero\oneone\onethree\onefour\onefive%
\twozero\twoone\twotwo\twofive\twoseven\twoeight\twonine%
\threezero\threeone\threetwo,%
345789\onetwo\onesix\oneseven\oneeight\onenine\twothree\twofour\twosix)
$$
In other words, the second row of the tabloid
records positions of occurrences of letter 1 in $p$.
Corresponding to the second row of the tabloid,
the RC-bijection proceeds as follows:
\begin{center}
\unitlength 12pt
\begin{picture}(31,2.5)(0,-0.5)
\put(0,1){$\emptyset$}
\put(1,1){$\xrightarrow{3}$}
\put(2.5,0.8){\yng(1)}
\put(4,1){1}
\put(5,1){$\xrightarrow{4}$}
\put(6.5,0.8){\yng(2)}
\put(9,1){0}
\put(10,1){$\xrightarrow{5}$}
\put(11.6,0.8){\yng(3)}
\put(15.3,1){$-1$}
\put(17,1){$\xrightarrow{7}$}
\put(18.5,-0.3){\yng(3,1)}
\put(22,1){$-1$}
\put(19.9,-0.2){3}
\put(23.7,1){$\xrightarrow{8}$}
\put(25,-0.3){\yng(4,1)}
\put(29.8,1){$-2$}
\put(26.4,-0.2){3}
\end{picture}
\end{center}
\begin{center}
\unitlength 12pt
\begin{picture}(31,4.5)(0,-0.5)
\put(0,3){$\xrightarrow{9}$}
\put(1.5,1.7){\yng(5,1)}
\put(7.2,3){$-3$}
\put(2.8,1.8){3}
\put(9,3){$\xrightarrow{12}$}
\put(10.7,0.6){\yng(5,1,1)}
\put(16.5,3){$-3$}
\put(12.1,1.9){6}
\put(12.1,0.8){3}
\put(18.2,3){$\xrightarrow{16}$}
\put(20,-0.5){\yng(5,1,1,1)}
\put(25.8,3){$-3$}
\put(21.4,1.9){8}
\put(21.4,0.8){6}
\put(21.4,-0.3){3}
\end{picture}
\end{center}
\begin{center}
\unitlength 12pt
\begin{picture}(31,4.5)(0,-0.5)
\put(0,3){$\xrightarrow{17}$}
\put(1.7,-0.5){\yng(5,2,1,1)}
\put(7.5,3){$-3$}
\put(4.3,2){5}
\put(3.1,0.8){6}
\put(3.1,-0.3){3}
\put(9.2,3){$\xrightarrow{18}$}
\put(10.7,-0.5){\yng(5,3,1,1)}
\put(16.5,3){$-3$}
\put(14.3,2){2}
\put(12.1,0.8){6}
\put(12.1,-0.3){3}
\put(18.2,3){$\xrightarrow{19}$}
\put(20,-0.5){\yng(5,4,1,1)}
\put(25.8,3){$-3$}
\put(24.6,1.9){$-1$}
\put(21.4,0.8){6}
\put(21.4,-0.3){3}
\end{picture}
\end{center}
\begin{center}
\unitlength 12pt
\begin{picture}(31,6.5)(0,-1.5)
\put(0,4){$\xrightarrow{23}$}
\put(1.7,-0.7){\yng(5,4,1,1,1)}
\put(7.5,4){$-3$}
\put(6.4,2.9){$-1$}
\put(3.1,1.8){13}
\put(3.1,0.7){6}
\put(3.1,-0.5){3}
\put(9.2,4){$\xrightarrow{24}$}
\put(10.7,-0.7){\yng(5,4,2,1,1)}
\put(16.5,4){$-3$}
\put(15.3,2.9){$-1$}
\put(13.2,1.8){8}
\put(12.1,0.6){6}
\put(12.1,-0.5){3}
\put(18.2,4){$\xrightarrow{26}$}
\put(20,-1.7){\yng(5,4,2,1,1,1)}
\put(25.8,4.1){$-3$}
\put(24.6,3){$-1$}
\put(22.5,1.9){8}
\put(21.4,0.8){14}
\put(21.4,-0.3){6}
\put(21.4,-1.4){3}
\end{picture}
\end{center}
In the above diagrams,
we attach each rigging on the right of
the corresponding row of the partition.
The last rigged configuration in the above sequence
gives the output of RC-bijection applied to $p$.
\hfill$\square$
\end{example}

\subsection{Kuniba--Takagi--Takenouchi's construction}
\subsubsection{Crystals $B_l$}
Let $B_l$ be the classical crystal \cite{Kas}
corresponding to the $l$-fold symmetric tensor
representation of $U_q'(A^{(1)}_1)$.
As the set, $B_l=\{(x_0,x_1)\in\mathbb{Z}_{\geq 0}^2\, |\,
x_0+x_1=l\}$.
Note that in usual notation, say that used in \cite{KTT},
our 0 and 1 are denoted by 1 and 2, respectively. 
We sometimes represent elements of $B_l$ by
semi-standard Young tableaux (without frame for simplicity).
More precisely, the element $(x_0,x_1)$ is also represented by
$\overbrace{0\cdots 0}^{x_0}
\overbrace{1\cdots 1}^{x_1}$.
In this notation we
usually denote $(1,0)$ and $(0,1)$
simply by 0 and 1, respectively.
We can define tensor product $B_l\otimes B_k$.
As the set, it is the product of the two sets
and it is known that we can define algebraic structure on them.
Then we can define the map
$R: B_l \otimes B_1 \rightarrow B_1 \otimes B_l$
which is compatible with the algebraic structure.
Explicitly, it is given by
\begin{align}
(x_0,x_1)\otimes 0 &\longmapsto 
\begin{cases}0 \otimes (l,0) &\hbox{if } (x_0,x_1)=(l,0)\\
1 \otimes (x_0+1,x_1-1) & \hbox{otherwise},
\end{cases}\\
(x_0,x_1)\otimes 1 &\longmapsto 
\begin{cases}1 \otimes (0,l) &\hbox{if } (x_0,x_1)=(0,l)\\
0 \otimes (x_0-1,x_1+1) & \hbox{otherwise}.
\end{cases}
\end{align}
$R$ is a bijection and called the combinatorial $R$ matrix.
We write the relation $R(u \otimes b)= b' \otimes u'$ 
simply as $u \otimes b \simeq b' \otimes u'$, and similarly for 
any consequent relation of the form
$a\otimes u \otimes b \otimes c \simeq
a \otimes b' \otimes u' \otimes c$.
When we consider elements of $p\in B_1^{\otimes L}$,
we sometimes omit symbols $\otimes$ and identify them
with paths.

\subsubsection{Time evolutions}
Given a path $p=b_1\otimes b_2\otimes\cdots\otimes b_L
\in B_1^{\otimes L}$, we define
the time evolution operators $T_l$ for $l\in\mathbb{Z}_{>0}$
of the PBBS in the following way.
Take $u_l:=(l,0)\in B_l$ and define $v_l\in B_l$ by
\begin{align}\label{def:v_l}
u_l\otimes p\simeq p_l^{\ast}\otimes v_l,\qquad
(p^\ast_l\in B_1^{\otimes L}).
\end{align}
In practice, we can compute the map $R$ step by step
as follows:
\begin{align}
u_l\otimes b_1\otimes b_2\otimes\cdots\otimes b_L\,
&\simeq
b_1^\ast\otimes u_l^{(1)}\otimes b_2\otimes\cdots\otimes b_L\simeq
b_1^\ast\otimes b_2^\ast\otimes u_l^{(2)}\otimes b_3
\otimes\cdots\otimes b_L\nonumber\\
&\simeq\cdots\simeq
b_1^\ast\otimes b_2^\ast\otimes\cdots\otimes b_L^\ast
\otimes u_l^{(L)},
\end{align}
where $u_l^{(i)}\in B_l$ and $u_l^{(L)}=v_l$.
This element $v_l$ depends on the path $p$.
Then, according to Section 2 of \cite{KTT}, we have
\begin{align}
v_l\otimes p\simeq
p^{\ast\ast}_l\otimes v_l.
\end{align}
Using this non-trivial relation, we define
\begin{align}
T_l(p):=p^{\ast\ast}_l.
\end{align}
Since combinatorial $R$ matrix acts trivially on $B_1\otimes B_1$,
$T_1$ is just the cyclic shift operator:
$T_1(b_1\otimes\cdots\otimes b_{L-1}\otimes b_L)=
b_L\otimes b_1\otimes\cdots\otimes b_{L-1}$.
The commutativity $T_lT_k=T_kT_l$ for all
$l,k\in\mathbb{Z}_{>0}$ holds due to the
Yang--Baxter relation (\cite{KTT} Theorem 2.2).
In general, there is sufficiently large $N$
so that $T_N=T_{N+1}=\cdots =:T_\infty$.
This definition of $T_\infty$ coincides with
the definition in Section \ref{sec:10elimination}
(Example 2.7 of \cite{KTT}).

\subsubsection{Action variable}
Consider the path $p=b_1\otimes\cdots\otimes b_L$
of weight $\lambda=(\lambda_1,\lambda_2)$, where
$\lambda_i=\#\{a\, |\,b_a=i\}$ satisfying $\lambda_1\geq\lambda_2$.
Then we can always find $d$ $(1\leq d\leq L)$ such that
$T_1^d(p)=b_1'\otimes\cdots\otimes b_L'$
has the following property:
for any $1\leq i\leq L$,
\begin{align}\label{Yamanouchi}
\#\{a\, |\,b_a'=1,a\in [i,L]\}\leq
\#\{a\, |\,b_a'=0,a\in [i,L]\}.
\end{align}
Let the rigged configuration corresponding to this $T_1^d(p)$
be $(\nu ,J)$.
We define the {\it action variable} of the path $p$ by $\nu$.
{}From the condition (\ref{Yamanouchi}),
we see that the configuration part of the RC data
corresponding to $T_1^d(p)^{\otimes N}$ is
$\overbrace{\nu\cup\nu\cup\cdots\cup\nu}^N$.
Due to Corollary 3.5 of \cite{KTT},
$\nu$ is conserved under arbitrary time evolution.
In the following,
for the given action variable $\nu$,  we put
\begin{align}
H=\{i\in\mathbb{Z}_{\geq 1}\, |\,m_i\neq 0\}
=\{i_1<i_2<\cdots <i_s\}.
\end{align}

\subsubsection{Angle variable}
Let us define angle variable according to \cite{KTT}.
In order to explain motivation of the construction, let us
mention the following simple property.
\begin{lemma}
Let $p$ be a path of length $L$ satisfying
the condition (\ref{Yamanouchi}), and $(\nu ,J_{\alpha ,i})$
be the corresponding rigged configuration.
Then the rigged configuration corresponding to the path
$\pi_N =\overbrace{p\otimes p\otimes\cdots\otimes p}^N$
is equal to
\begin{align}
\nu(\pi_N)&=\nu\cup\nu\cup\cdots\cup\nu,\\
J(\pi_N)&=\coprod_{a=1}^N\{J_{k+am_i,i}\, |\,
i\in H, 1\leq k\leq m_i\},
\end{align}
where $J_{k+am_i,i}:=J_{k,i}+aP_i(\nu)$.
\hfill$\square$
\end{lemma}

Motivated by this lemma, define $\bar{\mathcal{J}}$
as follows:
\begin{align}
\bar{\mathcal{J}}&=\bar{\mathcal{J}}(\nu )=
\mathcal{J}_{i_1}\times\mathcal{J}_{i_2}
\times\cdots\times
\mathcal{J}_{i_s},\\
\mathcal{J}_{i}&=
\{(J_{k,i})_{k\in\mathbb{Z}}\, |\,
J_{k,i}\in\mathbb{Z},
J_{k,i}\leq J_{k+1,i},
J_{k+m_i,i}=J_{k,i}+P_i(\nu)
\mbox{ for all }k\} .
\end{align}
\begin{definition}
For $j\in\mathbb{Z}_{\geq 1}$, define the map
$\sigma_j :\mathcal{J}_i\longrightarrow\mathcal{J}_i$ by
\begin{align}
\sigma_j:
(J_{k,i})_{k\in\mathbb{Z}}\longmapsto
(J_{k,i}')_{k\in\mathbb{Z}},\qquad
J_{k,i}'=J_{k+\delta_{i,j},i}+2\min (i,j).
\end{align}
We extend $\sigma_j$ to the map
$\bar{\mathcal{J}}\longrightarrow\bar{\mathcal{J}}$
by $\sigma_j(\bar{\mathcal{J}})=
\sigma_j(\bar{\mathcal{J}}_{i_1})\times
\sigma_j(\bar{\mathcal{J}}_{i_2})\times
\cdots\times
\sigma_j(\bar{\mathcal{J}}_{i_s})$.
\hfill$\square$
\end{definition}
The operators $\sigma_j$ are the key to define angle variable.
By definition, we have $\sigma_j\sigma_k=\sigma_k\sigma_j$
for any $j,k\in\mathbb{Z}_{\geq 1}$.
Therefore the set $\mathcal{A}:=\{\sigma_{i_1}^{n_1}
\sigma_{i_2}^{n_2}\cdots
\sigma_{i_s}^{n_s}\, |\,
n_1,n_2,\cdots,n_s\in\mathbb{Z}\}$
is abelian group.
\begin{definition}
$(1)$ Define equivalence relation $\simeq$ for elements of
$\bar{\mathcal{J}}$ by
$J\simeq K$ $(J,K\in\bar{\mathcal{J}})$
if there exist $\sigma\in\mathcal{A}$
such that $J=\sigma(K)$.\\
$(2)$ Define the set $\mathcal{J}=\mathcal{J}(\nu)$ by
\begin{align}
\mathcal{J}=\bar{\mathcal{J}}/\simeq.
\end{align}
An element of $\mathcal{J}$ is called angle variable.
\hfill$\square$
\end{definition}

Now we shall explain how to compute the operators
$\sigma_k$ based on explicit example.
To begin with, we have to embed a sequence of the riggings
$(J_{k,i})_{1\leq k\leq m_i}$ attached to length $i$
rows of $\nu$ into the set $\mathcal{J}_i$
by the map
\begin{align}
\iota :(J_{k,i})_{1\leq k\leq m_i}
\longrightarrow (J_{k,i})_{k\in\mathbb{Z}}
\end{align}
using periodicity $J_{k+m_i,i}=J_{k,i}+P_{i}(\nu)$.
\begin{example}
Consider the following length $L=19$ paths:
\begin{align}
q=0011010001110100011,\qquad
q'=0001110100011001101.
\end{align}
Notice that $q=T_1^6(q')$.
The rigged configurations corresponding to
$q$ and $q'$ are
\begin{center}
\unitlength 12pt
\begin{picture}(25,5.5)
\put(0,2.2){$q\longmapsto (\mu ,J)=$}
\put(6.5,0){$\Yboxdim 12pt \yng(3,2,2,1,1)$}
\put(7.8,0.1){2}
\put(7.8,1.1){6}
\put(8.8,2.1){0}
\put(8.8,3.1){3}
\put(9.8,4.1){0}
\put(10.7,2){,}
\put(13,2.2){$q'\longmapsto (\mu ,J')=$}
\put(20,0){$\Yboxdim 12pt \yng(3,2,2,1,1)$}
\put(21.3,0.1){4}
\put(21.3,1.1){9}
\put(22.3,2.1){3}
\put(22.3,3.1){3}
\put(23.3,4.1){0}
\put(24,2.2){.}
\end{picture}
\end{center}
We compute actions of $\sigma_2$ and $\sigma_1$
on the rigged configuration $(\mu ,J)$ corresponding
to $q$.
Since $\mu =(3,2,2,1,1)$, we have
$Q_3(\mu)=9$, $Q_2(\mu)=8$ and $Q_1(\mu)=5$.
Therefore the vacancy numbers are
$P_3(\mu)=1$, $P_2(\mu)=3$ and $P_1(\mu)=9$.
We would like to explain the following computations:
\begin{center}
\unitlength 12pt
\begin{picture}(30,5.5)
\put(0,0){$\Yboxdim 12pt \yng(3,2,2,1,1)$}
\put(1.3,0.1){2}
\put(1.3,1.1){6}
\put(2.3,2.1){0}
\put(2.3,3.1){3}
\put(3.3,4.1){0}
\put(5,2.2){$\stackrel{\sigma_2}{\longmapsto}$}
\put(8,0){$\Yboxdim 12pt \yng(3,2,2,1,1)$}
\put(9.3,0.1){4}
\put(9.3,1.1){8}
\put(10.3,2.1){7}
\put(10.3,3.1){7}
\put(11.3,4.1){4}
\put(13,2.2){$\stackrel{\sigma_1}{\longmapsto}$}
\put(16,0){$\Yboxdim 12pt \yng(3,2,2,1,1)$}
\put(17.3,0.1){10}
\put(17.3,1.1){15}
\put(18.3,2.1){9}
\put(18.3,3.1){9}
\put(19.3,4.1){6}
\put(21,2.2){$=\quad 6\,\,\,\, +$}
\put(26,0){$\Yboxdim 12pt \yng(3,2,2,1,1)$}
\put(27.3,0.1){4}
\put(27.3,1.1){9}
\put(28.3,2.1){3}
\put(28.3,3.1){3}
\put(29.3,4.1){0}
\end{picture}
\end{center}

$\bullet$ Action of $\sigma_2$ is divided into the following
two parts.
(a) For the riggings associated to length 2 rows, i.e., $(0,3)$
we embed them into $\mathcal{J}_2$ by $\iota$ as follows:
\begin{align}
(0,3)\longmapsto
(\cdots ,-6,-3,-3,0,\underline{0,3,}3,6,6,9,\cdots).
\end{align}
Here the original $(0,3)$ is regarded as $(J_{1,2},J_{2,2})$
which is enlarged with periodicity $P_2(\mu)=3$.
As a result of $\sigma_2$, we take riggings
$(J_{2,2},J_{3,2})=(3,3)$.
(b) For the riggings associated to length $i$ rows,
we further add $2\min (i,2)$.

$\bullet$ Similarly, action of $\sigma_1$ is
computed as follows.
(a) For the riggings associated to length 1 rows, i.e., $(4,8)$
we embed them into $\mathcal{J}_1$ by $\iota$ as follows:
\begin{align}
(4,8)\longmapsto
(\cdots ,-14,-10,-5,-1,\underline{4,8,}13,17,22,26,\cdots).
\end{align}
Here the original $(4,8)$ is regarded as $(J_{1,1},J_{2,1})$
which is enlarged with periodicity $P_1(\mu)=9$.
As a result of $\sigma_1$, we take riggings
$(J_{2,1},J_{3,1})=(8,13)$.
(b) For the riggings associated to length $i$ rows,
we further add $2\min (i,1)$.

To summarize, we have $\sigma_1\sigma_2 (J)=6+J'$.
This computation has the following interpretation \cite{S2}.
(a) 6 corresponds to the exponent in $T_1^6(q')=q$.
(b) As for $\sigma_1\sigma_2$, we see that
in order to get $q$ from $q'$,
we have to move six letters $001101$ form right of $q'$ to left.
These six letters contains two solitons of length 2 and 1,
which explains the form of composition $\sigma_1\sigma_2$.
\hfill$\square$
\end{example}

\subsubsection{Direct and inverse scattering transform}
Recall that we can choose an integer $d$ such that
$p'=T_1^{d}(p)$ satisfies condition (\ref{Yamanouchi}).
We define direct scattering transform $\Phi$ by
\begin{equation}
\begin{split}
\Phi :\quad&\mathcal{P}(\nu)\longrightarrow
\mathbb{Z}\times\tilde{\mathcal{P}}(\nu)\longrightarrow
\,\,\,\,\bar{\mathcal{J}}(\nu)\,\,\,\longrightarrow
\,\,\,\,\,\mathcal{J}(\nu)\\
&\,\,\,\, p\quad \longmapsto
\,\,\,\, (d,p')\,\,\,\,\longmapsto
\iota (J)+d\longmapsto
[\iota (J)+d]
\end{split}
\end{equation}
where $\tilde{\mathcal{P}}(\nu )$ is the subset of
$\mathcal{P}(\nu )$ consisting of all paths
satisfying condition (\ref{Yamanouchi}).
In \cite{KTT}, it was shown that the inverse map
$\Phi^{-1}$ is well-defined.
\begin{theorem}[\cite{KTT}, Theorem 3.12]
Define $T_l$ on $\mathcal{J}(\nu )$ by
\begin{align}
T_l:(J_{k,i})_{k\in\mathbb{Z},i\in H}
\longmapsto
(J_{k,i}+\min (i,l))_{k\in\mathbb{Z},i\in H}.
\end{align}
Then the following diagram commutes:
\begin{equation}
\begin{CD}
{\mathcal P}(\nu) @>{\Phi}>> {\mathcal J}(\nu) \\
@V{T_l}VV @VV{T_l}V\\
{\mathcal P}(\nu) @>{\Phi}>> {\mathcal J}(\nu) 
\end{CD}
\end{equation}
\hfill$\square$
\end{theorem}

\section{Interpretation of 10-elimination in terms
of the rigged configurations}\label{sec:main}
Let $p$ be a path of type $B_1^{\otimes L}$.

\begin{theorem}\label{th:main}
Let $x_i^{(j)}$ be a position of the $i$-th 0-soliton
located on the $E^{L_j}(p)$.
Let $T(p)$ be a two-row tabloid corresponding to the path $p$,
and $J_1$, $\cdots$, $J_g$ be the riggings corresponding to tabloid
$T(p)$ under the RC-bijection.
Then,
\begin{align}\label{eq:main}
x_i^{(j)}=J_{i,L_j}+L_j.
\end{align}
\end{theorem}
{\bf Proof.}
We use induction on the length of path.
Let the configuration corresponding to length
$L$ path under the RC-map be $\nu_L$.
The case $L=1$ is clear.
Let $L>1$.
Let $\pi$ be a path of length $L+1$.
Then one can distinguish two cases:
$\pi =p0$ or $\pi = p1$.
In the first case, the positions of 0-solitons as well as
RC data remain the same.

In the case $\pi =p1$ we have to distinguish further two cases:
$\pi =p'01$ and $\pi =p''11$.
In the case $\pi =p'01$, the rigged configuration
corresponding to the path $\pi$ can be obtained from that
corresponding to the path $\pi'=p'0$ of length $L$
by creating the new length 1 singular string:
The rigging corresponding to the new singular string
is equal to
\begin{align}
J_{m_1 (\nu_{L+1}),1}=
L+1-2(\ell (\nu_L)+1).
\end{align}
Therefore
\begin{align}
J_{m_1 (\nu_{L+1}),1}+1=L-2\ell (\nu_L).
\end{align}
But it is easily seen that exactly one new 0-soliton
appears in the position $L-2\ell (\nu_L)$.
In the case $\pi=p'11$ we observe that the position of
all 0-soliton do not change.
Now we have to check that all numbers
$J_{i,L_j}+L_j$ do not change as well.
Indeed, according to the RC-algorithm,
we have to add one box to the singular string
in the highest position.
After creating the new singular string, its rigging
becomes equal to the old one minus 1.
But the length of the new string becomes
bigger by exactly 1.
That means that the modified riggings
(i.e., sum of length of row and the corresponding rigging)
do not change.
\hfill$\blacksquare$

\bigskip

Let $X=\{x_0<x_1<\cdots <x_N<\cdots\}$ be an ordered set.
For any $x\in X,$ $x=x_i$, define $x^{\pm}=x_{i\pm 1}\in X,$ if $i \ge 1,$
and $x_{0}^{-}=x_{0}.$ Let $\alpha$ be a composition.
A standard X-tabloid of the 
shape $\alpha$ is a filling of the diagram
corresponding to the composition $\alpha$ 
by the elements from the totally ordered set $X \setminus\{x_0\}$ such that 
the elements along each row are strictly increasing, and the all elements of 
the filling in question are distinct.
An element $x \in T$ is called descent (resp. acsent) if either $x \not= x_1$ 
and the element $x^+$ belongs to $T$, and located in the 
tabloid $T$ below (resp. above) the element $x,$ ~~or $x=x_1$ and it
does't belong to the first row of the tabloid $T.$
We denote by $Des(T)=\{x^{+} |\, x \in T ~\mbox{is a descent} \}$
(resp. $Acs(T)=\{x^{-} |\,x \in T ~\mbox{is an acsent}\}$).
Let $d(T)=\# \, |\,Des(T)\, |\,$ (resp. $a(T)=\# \, |\,Acs(T)\, |\,$).

To continue, let us recall that under the rigged configuration bijection to 
a (semi) standard tabloid $T$ of shape $\alpha=(\alpha_1,\ldots,\alpha_r)$ one 
associates a collection of partitions $(\nu^{(1)},\ldots,\nu^{(r-1)})$ such 
that $| \nu^{(k)}|=\sum_{ j \ge k+1} \alpha_j,$~$k=1,\ldots,r-1,$ and a set of
non-negative integers (riggings) $\{J_{i,j}^{(k)} \},$ ~
$j=\nu_1^{(k)}, \nu_2^{(k)},\ldots,$~~~$1 \le i \le \#|s, j=\nu_s^{(k)}|,$ 
with certain constraints, see e.g. \cite{KR} for details. Let $T$ be a (semi)
standard tabloid, we will call the partition $\nu^{(1)}$
which appears after applying the RC-bijection to $T$, to be the first 
configuration corresponding to the tabloid $T.$

Our next goal is to define recursively a sequence of
standard tabloids $T_0=T$, $T_1,T_2,\cdots$
to be used to define the shape of the first configuration
corresponding to a tabloid $T$ under
the RC-bijection.
Thus  let us describe how to obtain tabloid
$T_i$ from that $T_{i-1}$, $i\geq 1$.
Namely, consider the descent set $Des(T_{i-1})$,
and denote by $T_i$ the tabloid obtained from that $T_{i-1}$
by deleting all entries of $T_{i-1}$ which are
equal to either $x$ or $x^-$ for all
$x\in Des (T_{i-1})$.
\begin{theorem}[\cite{K1}]
Let $T$ be a standard tabloid,
denote by $\nu$ the shape of the first configuration
corresponding to tabloid $T$ under the RC-bijection.
Then
\begin{align}
\nu_i'=\# \, |\,Des (T_{i-1})\, |\,,\qquad
i=1,2,\cdots ,
\end{align}
where $\nu'$ is the transposition of $\nu$.
\hfill$\square$
\end{theorem}
In the case when a tabloid $T$ has two rows and corresponds to a positive weight path,
one can use a similar construction by using the sets
$Acs(T_{i-1})$ instead of those $Des(T_{i-1})$.
It is easy to see that in the $A^{(1)}_1$ case, for positive weight paths 
we also have
\begin{align}
\nu_i'=\# \, |\,Acs(T_{i-1})\, |\,,\qquad
i\geq 1.
\end{align}

As a corollary, we see that in the case under consideration, the 10-elimination 
algorithm produces the same shape as the
RC-algorithm does:
\begin{align}
\nu =(L_1^{m_{L_1}},L_2^{m_{L_2}},\cdots,L_s^{m_{L_s}}),
\end{align}
where $\nu$ is the configuration associated with
the given path under RC-algorithm,
and $L_i^{m_{L_i}}$ $(i=1,2,\cdots,s)$ are
determined by the 10-elimination 
algorithm.

To summarize, we obtain the following main theorem of our paper.
\begin{theorem}\label{th:main}
Let $p$ be a $B_1^{\otimes L}$ path of $A^{(1)}_1$.
Then we have
\begin{align}
(\nu ,J)=\{(L_i,x_\alpha^{(i)}-L_i)
\}_{1\leq i\leq s,1\leq \alpha\leq m_{L_i}}
\end{align}
where $(\nu ,J)$ is the rigged configuration associated with
the given path under RC-algorithm,
and $L_i^{m_{L_i}}$, $x_\alpha^{(i)}$ are the parameters
determined by the 10-elimination 
algorithm.
\hfill$\square$
\end{theorem}

\begin{example}
Consider the length 32 path $p$
treated in Example \ref{ex:MIT}.
Let us compare (\ref{MITdata1}) and (\ref{MITdata2})
with the result of Example \ref{ex:MIT2}.
Indeed, (\ref{MITdata1}) coincides with the partition
$(5,4,2,1,1,1)$ obtained in Example \ref{ex:MIT2}.
As for (\ref{eq:main}), we compute
\begin{align*}
\{x^{(1)}_1-L_1, x^{(2)}_1-L_2, x^{(3)}_1-L_3,
x^{(4)}_1-L_4, x^{(4)}_2-L_4, x^{(4)}_3-L_4\}
=\{-3,-1,8,14,6,3\}
\end{align*}
which coincides with the riggings obtained in
Example \ref{ex:MIT2}.

We now demonstrate how the rigged configuration
essentially obtained by 10-elimination
will yield the solution of the initial value problem.
As an example, we compute $T_\infty^{10000}(p)$.
Under $T_\infty^{10000}$, the riggings behave as
\begin{center}
\unitlength 12pt
\begin{picture}(21,6.5)
\put(0,2.8){$T_\infty^{10000}:$}
\put(4,0){$\Yboxdim 12pt \yng(5,4,2,1,1,1)$}
\put(9.2,5.1){$-3$}
\put(8.2,4.1){$-1$}
\put(6.2,3.1){8}
\put(5.2,2.1){14}
\put(5.2,1.1){6}
\put(5.2,0.1){3}
\put(12,2.8){$\longmapsto$}
\put(15,0){$\Yboxdim 12pt \yng(5,4,2,1,1,1)$}
\put(20.2,5.1){49997}
\put(19.2,4.1){39999}
\put(17.2,3.1){20008}
\put(16.2,2.1){10014}
\put(16.2,1.1){10006}
\put(16.2,0.1){10003}
\end{picture}
\end{center}
To save space, we hereafter abbreviate such computation as
\begin{align}\label{eq:ex:1}
(-3,-1,8,14,6,3)
\xrightarrow{T_\infty^{10000}}
(49997,39999,20008,10014,10006,10003).
\end{align}
The remaining task is
to find a proper representative
in the equivalence class corresponding to the
rigging vector on the right hand side so that
we can apply the rigged configuration bijection.
To begin with, recall that the action of
the cyclic shift operator $T_1^d$ on the rigging vector
simply means constant shift $d=(d,d,\cdots,d)$
mod $L$, where $L$ is the length of the system.
Then we can always choose a suitable element
$\sigma\in\mathcal{A}$ and a constant vector
$d$ so that the riggings
$J'=\sigma (J)-d$ satisfy
\begin{align}
0\leq J_{1,i}'\leq J_{2,i}'\leq\cdots\leq
J_{m_i,i}'\leq P_i(\nu),
\qquad\forall i\in H.
\end{align}
Let us explain the method based on the example
under consideration.
In this case, we have $L=32$, $\nu=(5,4,2,1,1,1)$,
and thus $P_5(\nu)=4$, $P_4(\nu)=6$,
$P_2(\nu)=14$ and $P_1(\nu)=20$.
Remind that by acting $\sigma_k$, the riggings
corresponding to rows with length strictly bigger
than $k$ gains uniform shift by $2\cdot k$.
{}From this property, we can adjust the riggings
successively from longer rows to shorter ones.
In the above example, suppose that we apply $\sigma_4^N$
for some integer $N$.
Then the rigging corresponding to length 5 row
gains $N\times (2\cdot 4)=8N$
whereas the rigging corresponding to length 4
row gains $N\times (P_4(\nu)+2\cdot 4)
=14N$ (remind that $m_4(\nu)=1$).
By choosing $N=1667$, we get
\begin{align}\label{eq:ex:2}
\mbox{RHS of }
(\ref{eq:ex:1})
\xrightarrow{\sigma_4^{1667}}
(63333,63337,26676,13348,13340,13337).
\end{align}
Similarly, from $P_2(\nu)=14$, we get
\begin{align}\label{eq:ex:3}
\mbox{RHS of }
(\ref{eq:ex:2})
\xrightarrow{\sigma_4^{2619}}
(73809,73813,73818,18586,18578,18575).
\end{align}
Finally, from $P_1(\nu)=20$ and $m_1(\nu)=3$,
we get
\begin{align*}
\mbox{RHS of }
(\ref{eq:ex:3})
\xrightarrow{\sigma_1^{3\cdot 2762}}\,&
(90381,90385,90390,90398,90390,90387)\\
=\,&\, 90381+(0,4,9,17,9,6).
\end{align*}
Such adjustment is always possible because
the vacancy numbers corresponding to
all rows of $\nu$ except for the longest rows
are strictly positive
(Lemma 3.9 of \cite{KTT}).
On the rigged configuration $(\nu ,J')$,
$J'=(0,4,9,17,9,6)$, we can apply
the RC-bijection (see Appendix A to \cite{KTT}),
or equivalently the explicit formula (\ref{eq:xk})
to get the corresponding path $p'$:
$$p'=00000111011100100001100011101100.$$
Since $90381\equiv 13\pmod{32}$, we have
\begin{align*}
T_1^{90381}(p')=T_1^{13}(p')=
11000111011000000011101110010000,
\end{align*}
which reconstructs the computation
in the final part of \cite{MIT}.
\hfill$\square$
\end{example}
\begin{remark}
There are explicit formula for the initial
value problem of the PBBS in terms of the
ultradiscrete (or tropical) Riemann theta
function \cite{KS2,KS3}.
The formulae are direct consequence of
the formula (\ref{eq:xk}) bellow
and are logically independent to
the main claims of \cite{KTT}.
These formulae are parametrized by
the rigged configurations and thus can be completely
determined by the data $L_j$ and $x_i^{(j)}$ derived from the
10-elimination procedure.
\hfill$\square$
\end{remark}
\begin{remark}
All considerations in this section can be used
for linear systems as well as the BBS.
Let $p$ be an arbitrary linear path.
Then we can embed $p$ into semi-infinite
system by
\begin{align}
p\hookrightarrow p\otimes 1\otimes1\otimes1\otimes\cdots .
\end{align}
Then we can use Theorem \ref{th:main} and
considerations following it without any changes.

In \cite{KSY,S0} complete solution of the
initial value problem of the BBS was derived for
all $\bigotimes_i B_{l_i}$ type paths of $A^{(1)}_n$.
We quote here the formula in the special case of
$B_1^{\otimes L}$ of $A^{(1)}_1$.
Let the path be $p=b_1\otimes b_2\otimes
\cdots\otimes b_L\in B_1^{\otimes L}$
and denote the element $b_k$ by $(1-x(k),x(k))$,
i.e., if $b_k =1$ then $x(k)=1$ and $x(k)=0$ otherwise.
Denote the corresponding rigged configuration
by $(\nu ,J)=(\nu_i,J_i)_{i=1}^g$.
Here we do not assume that all $\nu_i,$~$i=1,2,\ldots$  are distinct.
Then
\begin{align}
x(k) &= \tau_0(k)-\tau_0(k-1)-\tau_1(k)+\tau_1(k-1),
\label{eq:xk}\\
\tau_r(k) &= -\min_{n \,\in \{0,1\}^g}
\left\{\sum_{i=1}^g(J_i+r\nu_i -k)n_i 
+ \sum_{i,j =1}^g\min(\nu_i,\nu_j)n_in_j\right\}
\quad (r=0,1),\label{eq:tau}
\end{align}
where $n=(n_1,n_2,\cdots,n_g)$.
As for the time evolution,
the rigged configuration corresponding to
$T_l(p\otimes 1\otimes1\otimes\cdots)$ is
$(\nu_i,J_i+\min (\nu_i,l))_{i=1}^g$
and we can plug this into (\ref{eq:xk}) directly.
Note that these formulae are parametrized by
the rigged configurations.
Therefore these formulae are completely determined by
the data $L_j$ and $x_i^{(j)}$ derived from the
10-elimination procedure.
\hfill$\square$
\end{remark}

\vspace{5mm}
\noindent
{\bf Acknowledgements:}
The work of RS is supported by
the Core Research for Evolutional Science and Technology
of Japan Science and Technology Agency.


\begin{thebibliography}{99}
\bibitem{Fuk} K.~Fukuda,
Box-ball systems and Robinson--Schensted--Knuth correspondence,
J. Algebraic Combin.   {\bf 19}  (2004) 67--89,
arXiv:math/0105226.

\bibitem{FOY} K.~Fukuda, M.~Okado and Y.~Yamada,
Energy functions in box-ball systems,
Int. J. Mod. Phys. {\bf A15} (2000) 1379--1392,
arXiv:math/9908116.

\bibitem{HHIKTT}
G.~Hatayama, K.~Hikami, R.~Inoue, A.~Kuniba,
T.~Takagi and T.~Tokihiro,
The $A^{(1)}_M$ automata related to crystals of symmetric tensors,
J. Math. Phys. {\bf 42} (2001) 274--308, arXiv:math/9912209.

\bibitem{IT1}
R.~Inoue and T.~Takenawa,
Tropical spectral curves and integrable cellular automata,
Int. Math. Res. Notices {\bf 2008} (2008)
Article ID rnn019 (27pp),
arXiv:0704.2471.

\bibitem{IT2}
R.~Inoue and T.~Takenawa,
A tropical analogue of Fay's trisecant identity
and the ultra-discrete periodic Toda lattice,
arXiv:0806.3318.

\bibitem{Kas} M.~Kashiwara,
On crystal bases of the $q$-analogue of universal enveloping algebras,
Duke Math. J. {\bf 63} (1991) 465--516.

\bibitem{K1} A.~N.~Kirillov, On some properties of the Robinson-Schensted 
correspondence, Proceedings of Hayashibara conference on special functions, (1990)
122--126.

\bibitem{KR} A.~N.~Kirillov and N.~Yu.~Reshetikhin, The Bethe ansatz and the 
combinatorics of Young tableaux,~~Zap. Nauch. Sem. LOMI {\bf 155} (1986), 65-115, 
translation in Journal of Soviet Math. {\bf 41} (1988) 925--955.

\bibitem{KS1}
A.~N.~Kirillov and R.~Sakamoto,
Paths and Kostka--Macdonald polynomials,
arXiv:0811.1085.

\bibitem{KS2} A.~Kuniba and R.~Sakamoto,
The Bethe ansatz in a periodic box-ball system and the ultradiscrete
Riemann theta function,
J. Stat. Mech. (2006) P09005, 1--12, arXiv:math/0606208.

\bibitem{KS3} A.~Kuniba and R.~Sakamoto,
Combinatorial Bethe ansatz and ultradiscrete
Riemann theta function with rational characteristics,
Lett. Math. Phys. {\bf 80} (2007) 199--209,
arXiv:nlin/0611046.

\bibitem{KSY} A.~Kuniba, R.~Sakamoto and Y.~Yamada,
Tau functions in combinatorial Bethe ansatz,
Nucl. Phys. {\bf B786} (2007) 207--266,
arXiv:math/0610505.

\bibitem{KTT}
A.~Kuniba, T.~Takagi, and A.~Takenouchi,
Bethe ansatz and inverse scattering transform 
in a periodic box-ball system, 
Nucl. Phys. {\bf B747} (2006) 354--397,
arXiv:math/0602481.

\bibitem{MIT}
J.~Mada, M.~Idzumi and T.~Tokihiro,
On the initial value problem of
a periodic box-ball system,
J. Phys. A: Math. Gen. {\bf 39} (2006) L617--L623,
arXiv:nlin/0608037.

\bibitem{S0} R.~Sakamoto, 
Crystal interpretation of Kerov--Kirillov--Reshetikhin bijection II.
Proof for $\mathfrak{sl}_n$ case,
J. Algebraic Combin. {\bf 27} (2008) 55--98,
arXiv:math/0601697.

\bibitem{S2}
R.~Sakamoto,
Finding rigged configurations from paths,
RIMS Kokyuroku Bessatsu {\bf B11} (2009) 1--17,
arXiv:0804.2511. 

\bibitem{T} D.~Takahashi,
{On some soliton systems defined by using
boxes and balls}, 
in ``Proceedings of
the International Symposium on Nonlinear Theory and
Its Applications'' (NOLTA '93),
(1993) 555--558.

\bibitem{TS} D.~Takahashi and J.~Satsuma,
{A soliton cellular automaton},
J. Phys. Soc. Jpn. {\bf 59} (1990) 3514--3519.

\bibitem{TTMS}
T.~Tokihiro, D.~Takahashi, J.~Matsukidaira
and J.~Satsuma, {}From soliton equations to integrable
cellular automata through a limiting procedure,
Phys. Rev. Lett. {\bf 76} (1996) 3247--3250.

\bibitem{TTS}
M.~Torii, D.~Takahashi and J.~Satsuma,
Combinatorial representation of
invariants of a soliton cellular automaton,
Physica {\bf D92} (1996) 209--220.

\bibitem{Y}
Y.~Yamada,
A birational representation of Weyl group,
combinatorial $R$-matrix and discrete Toda
equation, in ``Physics and Combinatorics 2000'',
eds. A.~N.~Kirillov and N.~Liskova,
World Scientific (2001) 305--319.

\bibitem{YYT} D.~Yoshihara, F.~Yura and T.~Tokihiro,
Fundamental cycle of a periodic box-ball system,
J. Phys. A: Math. Gen. {\bf 36} (2003) 99--121.
arXiv:nlin/0208042.

\bibitem{YT} F.~Yura and T.~Tokihiro,
On a periodic soliton cellular automaton,
J. Phys. A: Math. Gen. {\bf 35} (2002) 3787--3801,
arXiv:nlin/0112041.
\end{thebibliography}
\end{document}